\documentclass[12pt,twoside]{article}

\usepackage{amsfonts}
\usepackage{amssymb}
\usepackage{amsmath}
\usepackage{graphicx}

\newcommand{\ignore}[1]{}

\def\@begintheorem#1#2{\par\bgroup{\sc #1\ #2. }\it\ignorespaces}
\def\@opargbegintheorem#1#2#3{\par\bgroup{\sc #1\ #2\ (#3). } \it\ignorespaces}
\def\@endtheorem{\egroup}
\newtheorem{theorem}{Theorem}[section]
\newtheorem{corollary}[theorem]{Corollary}
\newtheorem{lemma}[theorem]{Lemma}
\newtheorem{example}[theorem]{Example}
\newtheorem{proposition}[theorem]{Proposition}
\newtheorem{definition}[theorem]{Definition}
\newcommand{\bt}[1]{\begin{theorem}\label{#1}}
\newcommand{\bc}[1]{\begin{corollary}\label{#1}}
\newcommand{\bl}[1]{\begin{lemma}\label{#1}}
\newcommand{\be}[1]{\begin{example}\label{#1}}
\newcommand{\bp}[1]{\begin{proposition}\label{#1}}
\newcommand{\ba}[1]{\begin{algorithm}\rm\label{#1}}
\newcommand{\bd}[1]{\begin{definition}\rm\label{#1}}
\newcommand{\bpr}{\noindent {\em Proof. }}
\newcommand{\et}{\end{theorem}}
\newcommand{\ec}{\end{corollary}}
\newcommand{\el}{\end{lemma}}
\newcommand{\ee}{\end{example}}
\newcommand{\ep}{\end{proposition}}
\newcommand{\ed}{\end{definition}}
\newcommand{\epr}{{\ \vbox{\hrule\hbox{%
\vrule height1.3ex\hskip0.8ex\vrule}\hrule}}\\\par}

\def\R{\mathbb{R}}
\def\Z{\mathbb{Z}}

\def \G {{\cal G}}

\begin{document}

\title{\bf Robust Integer Programming}
\author{Shmuel Onn}

\date{\small Technion - Israel Institute of Technology, Haifa, Israel
\\onn@ie.technion.ac.il}

\maketitle

\begin{abstract}
We provide a complexity classification of four variants of robust
integer programming when the underlying Graver basis is given.
We discuss applications to robust multicommodity flows and multidimensional
transportation, and describe an effective parametrization of robust integer programming.
\end{abstract}

\section{Introduction}

Robust discrete optimization involves making optimal decisions under uncertainty
in the problem data, see \cite{KY} for a detailed development of this framework.
In this note we study robust integer programming problems where the uncertainty
occurs in the cost of the feasible points. More precisely, we consider a set of feasible
integer points which satisfy a system of linear inequalities in standard form given by
\begin{equation}\label{Set}
X \ := \left\{x\in\Z^n\,,\ Ax=b\,,\ l\leq x\leq u\right\}\ ,
\end{equation}
where $A$ is an integer $m\times n$ matrix, $b\in\Z^m$, and $l,u\in\Z^n$.
The uncertain cost to be minimized belongs to a finite nonempty set $C\subset\Z^n$
of potential {\em cost} vectors. The {\em worst-case cost} of point $x\in X$
is $\max\{cx:c\in C\}$. The {\em robust integer programming problem} is then
to find $x\in X$ attaining minimum worst-case cost, that is, to solve
$$\min_{x\in X}\max_{c\in C}cx\ .$$
When $C$ is a singleton, the problem reduces to standard linear
integer programming, and hence is no easier than the latter, which is NP-hard.
We therefore make some reasonable assumptions on $C$ and $X$ and study the complexity
of the problem under these assumptions. First, we assume that $C$ is either an
explicit {\em list} $C=\{c^1,\dots,c^k\}$ of given cost vectors,
or  a {\em box} $C=\{c\in\Z^n:d\leq c\leq e\}$ for some $d,e\in\Z^n$,
so that the cost of each decision variable lies in a given interval. Note that a
box $C$ has exponentially many elements and cannot be explicitly listed efficiently.
Second, we assume that the so-called {\em Graver basis} $\G(A)$ of the matrix $A$
defining $X$ is given; this object, defined in Section 2, plays a central role
in a recent theory of integer programming developed in \cite{Onn} and the
references therein. As illustrated in Section 3, the Graver basis can be computed
in polynomial time from $A$ for a variety of applications including multicommodity
flows and multidimensional transportation. Moreover, there is a parametrization
of {\em all integer programs} \cite{DO} under which the Graver basis can be
computed in polynomial time for each fixed parameter value.

The standard interpretation of robust discrete optimization is that the ``decision maker"
first chooses a point from $X$ and that ``nature" then chooses a cost from $C$.
However, the problem can more symmetrically be regarded as a two party game where
the $X$ player must pay $cx$ to the $C$ player. It is then also natural to consider
the ``dual" problem where the $C$ player makes its choice first, that is, the problem
$$\max_{c\in C}\min_{x\in X}cx\ .$$

We show the following complexity classification of robust integer programming.

\bt{Main}
The following complexity table holds for the robust integer programming problems
with the Graver basis of the matrix of $X$ given as part of the input:
\vskip.2cm
\begin{center}
\begin{tabular}{|c|cc|}
  \hline
   & list $C$ & box $C$ \\
  \hline
  $\min_X\max_C$ & NP-hard & polynomial time \\
  $\max_C\min_X$ & polynomial time & NP-hard \\
  \hline
\end{tabular}
\end{center}
\et

\vskip.2cm
We remark that when {\em costs} are replaced by {\em profits},
the corresponding problems reduce to the corresponding problems
with $C$ replaced by $-C:=\{-c:c\in C\}$ via the following relations,
resulting in the same complexities as in Theorem \ref{Main}:
$$\max_{x\in X}\min_{c\in C}cx\ =\ -\min_{x\in X}\max_{c\in -C}cx\ ,\quad
\quad\min_{c\in C}\max_{x\in X}cx\ =\ -\max_{c\in -C}\min_{x\in X}cx\ .$$

\vskip.2cm
In Section 2 we define the Graver basis and prove Theorem \ref{Main}. In Section 3
we discuss applications to robust multicommodity flows and robust multidimensional
transportation, and describe the parametrization of robust integer programming.

\section{Proof}

We begin by defining the Graver basis. Introduce a partial order $\sqsubseteq$ on $\R^n$
by $x\sqsubseteq y$ if $x_iy_i\geq 0$ and $|x_i|\leq |y_i|$ for $i=1,\ldots,n$.
The {\em Graver basis} of the integer $m\times n$ matrix $A$ is the finite set
$\G(A)\subset\Z^n$ of $\sqsubseteq$-minimal elements in $\{x\in\Z^n\,:\, Ax=0,\ x\neq 0\}$.

\vskip.2cm\noindent
{\em Proof of Theorem \ref{Main}.}
We prove the four entries of the complexity table one by one.
\vskip.2cm\noindent
{\em Entry (1,1):}
Consider the NP-complete problem of deciding, given $a\in\Z_+^n$, if there is a subset
$I\subseteq[n]:=\{1,\dots,n\}$ such that $\sum_{i\in I}a_i=\sum_{i\notin I}a_i$.
Let $a_0:=\sum_{i=1}^na_i$ and define $c^1,c^2\in\Z^{n+1}$ by $c^1:=(0,a)$ and
$c^2:=(a_0,-a)$. Next, define the set
$$X\ :=\ \{x\in\Z^{n+1}\ :\
{\bf 0}x=0,\ \ 1\leq x_0\leq 1,\ \ 0\leq x_i\leq 1,\ \ i=1,\dots,n\}\ ;$$
it is not hard to verify that the Graver basis of the matrix ${\bf 0}$
defining $X$ is given by $\G({\bf 0})=\pm\{{\bf 1}_i:0\leq i\leq n\}$
with ${\bf 1}_i$ the $i$th unit vector. Now, there is a bijection between $x\in X$
and $I\subseteq [n]$ with $I(x):=\{i:1\leq i\leq n,\ x_i=1\}$. For any $x\in X$,
$$\max\{c^1x,c^2x\}\ =\ \max\{\sum_{i\in I(x)}a_i,a_0-\sum_{i\in I(x)}a_i\}
\ =\ \max\{\sum_{i\in I(x)}a_i,\sum_{i\notin I(x)}a_i\}\ \geq \ {a_0\over 2}$$
with equality if and only $\sum_{i\in I(x)}a_i=\sum_{i\notin I(x)}a_i$.
So $\min_{x\in X}\max_{c\in C}cx={a_0\over2}$ if and only
if there is an $I$ with $\sum_{i\in I}a_i=\sum_{i\notin I(x)}a_i$,
and solution of the robust integer programming problem will enable
solution of the given NP-complete problem.
\vskip.2cm\noindent
{\em Entry (1,2):}
For of each feasible point $x\in X$ define $f(x):=\max_{c\in C}cx$. Then
$$f(x):=\max\{\sum_{i=1}^nc_ix_i: d_i\leq c_i\leq e_i\}
\ =\ \sum_{i=1}^n\max\{d_ix_i,e_ix_i\}\ =\ \sum_{i=1}^nf_i(x_i)\ ,$$
where $f_i(x_i):=\max\{d_ix_i,e_ix_i\}$. Since $f_i(x_i)$ is the maximum
of two univariate convex functions, it is also univariate convex.
Therefore $f(x)=\sum_{i=1}^nf_i(x_i)$ is separable convex.
So the robust integer programming problem is the integer program
$$\min_{x\in X}\max_{c\in C}cx\ =\
\min\{f(x)\,:\,x\in\Z^n,\ Ax=b,\ l\leq x\leq u\}$$
of minimizing a separable convex function $f$ over integer points satisfying lower and
upper bounds and a system of equations with defining matrix whose Graver basis $\G(A)$ is
given, which can be solved in polynomial time, see \cite[Theorem 3.12]{Onn} or \cite{HOW}.
\vskip.2cm\noindent
{\em Entry (2,1):}
For each cost vector $c\in C$ in the given list, consider the program
$$g(c)\ :=\ \min\{cx\,:\,x\in\Z^n,\ Ax=b,\ l\leq x\leq u\}\ ;$$
its objective function $cx$ is linear hence separable convex, and therefore, given $\G(A)$,
can be solved in polynomial time by \cite[Theorem 3.12]{Onn} again; then that $c\in C$
which attains maximum value $g(c)$ solves the given robust integer programming problem.
\vskip.2cm\noindent
{\em Entry (2,2):}
Consider again the NP-complete problem of deciding, given $a\in\Z_+^n$, if there is a
subset $I\subseteq[n]$ such that $\sum_{i\in I}a_i=\sum_{i\notin I}a_i$.
Let again $a_0:=\sum_{i=1}^n a_i$. Let
$$C\ :=\ \{c\in\Z^{n+2}\ :\
1\leq c_0\leq 1,\ \ 0\leq c_i\leq 1,\ \ i=1,\dots,n,\ \ 0\leq c_{n+1}\leq 0\}\ ,$$
Next let $m:=n+1$ and define an $m\times (n+2)$ matrix $A$ and vector $b\in\Z^m$ by
$$A\ :=\ \left(I_{n+1}\
      \begin{array}{r}
        a_0 \\
        -2a \\
      \end{array}
    \right),\quad
  b\ :=\ \left(
      \begin{array}{r}
        0 \\
        -a \\
      \end{array}
    \right)\ ,$$
and let
$$X\ :=\ \{x\in\Z^{n+2}\ :\ Ax=b,\ \
-|a_0|\leq x_i\leq |a_0|,\ \ i=0,\dots,n,\ \ 0\leq x_{n+1}\leq 1\}\ .$$
It is not hard to verify that the
Graver basis of $A$ is given by $\G(A)=\pm\{(-a_0,2a,1)\}$.

Now, note that the value of $x_{n+1}$ determines the value of the other $x_i$ via the
system $Ax=b$ so that $X=\{(0,-a,0),(-a_0,a,1)\}$. Next, there is a bijection between
$c\in C$ and $I\subseteq [n]$ with $I(c):=\{i:1\leq i\leq n,\ c_i=1\}$. For each $c\in C$,
$$\min_{x\in X}cx\ =\ \min\{-\sum_{i\in I(c)}a_i,-a_0+\sum_{i\in I(c)}a_i\}
\ =\ \min\{-\sum_{i\in I(c)}a_i,-\sum_{i\notin I(c)}a_i\}\ \leq \ -{a_0\over 2}$$
with equality if and only $\sum_{i\in I(c)}a_i=\sum_{i\notin I(c)}a_i$.
So $\max_{c\in C}\min_{x\in X}cx=-{a_0\over2}$ if and only
if there is an $I$ with $\sum_{i\in I}a_i=\sum_{i\notin I}a_i$,
and solution of the robust integer programming problem
will enable solution of the given NP-complete problem.
\epr

\section{Applications}

An {\em $(r,s)\times t$ bimatrix} is a matrix $A$ consisting of two blocks 
$A_1$, $A_2$, with $A_1$ its $r\times t$ submatrix consisting of the first $r$ 
rows and $A_2$ its $s\times t$ submatrix consisting of the last $s$ rows. 
The {\em $n$-fold product} of $A$ is the following $(r+ns)\times nt$ matrix,
$$A^{(n)}\quad:=\quad
\left(
\begin{array}{cccc}
  A_1    & A_1    & \cdots & A_1    \\
  A_2    & 0      & \cdots & 0      \\
  0      & A_2    & \cdots & 0      \\
  \vdots & \vdots & \ddots & \vdots \\
  0      & 0      & \cdots & A_2    \\
\end{array}
\right)\quad .
$$
For each fixed bimatrix $A$, the Graver basis $\G(A^{(n)})$ of the $n$-fold product 
of $A$ can be computed in time which is polynomial in $n$, see \cite[Theorem 4.4]{Onn} 
or \cite{DHOW}.
This has a variety of applications including multicommodity flows, multidimensional
transportation, and more generally any integer program via a suitable parametrization.
We now discuss consequences of this to the robust counterparts of these applications.

The (integer) {\em multicommodity flow problem} is as follows.
There are $l$ types of discrete commodities, $m$ suppliers, and $n$ consumers.
Supplier $i$ has supply $s^k_i$ in commodity $k$. Consumer $j$ has
demand $d^k_j$ in commodity $k$. Channel $(i,j)$ has capacity $u_{i,j}$
which is an upper bound on the total flow of all commodities on that
channel. There is a cost $c^k_{i,j}$ per unit flow of commodity $k$ on channel $(i,j)$.
The non-robust problem is to find a multicommodity flow consisting of $x^k_{i,j}$
units of flow of commodity $k$ on channel $(i,j)$ for all $i,j,k$, satisfying
the supply, demand, and capacity constraints, and attaining minimum cost
$\sum_{i,j,k}c^k_{i,j} x^k_{i,j}$, that is, solve
$$\min\{cx\ :\ x=(x^k_{i,j})\in\Z_+^{l\times m\times n}\,,\
\sum_j x^k_{i,j}=s^k_i\,,\ \ \sum_i x^k_{i,j}=d^k_j\,,\ \
\sum_{k=1}^l x^k_{i,j}\leq u_{i,j}\}\ .$$
The problem is NP-hard already for $l=2$ commodities or $m=3$ suppliers.
However, it is natural to have relatively small numbers of commodities and
suppliers but very large number of consumers, and we have
the following consequence of our theorem.
\bc{Flows}
For fixed $l$ commodities and $m$ suppliers, the robust multicommodity flow
problem $\min_X\max_C$ with box $C$ or $\max_C\min_X$ with list $C$ is
polytime solvable.
\ec
\bpr
Introduce a slack commodity $0$ with new variables
$x^0_{i,j}:=u_{i,j}-\sum_{k=1}^l x^k_{i,j}$ representing the slack flow
on each channel $(i,j)$, with cost $c^0_{i,j}:=0$, and suitable slack
supplies $s^0_i:=\sum_j u_{i,j}-\sum_{k=1}^ls^k_i$ and slack
demands $d^0_j:=\sum_i u_{i,j}-\sum_{k=1}^ld^k_j$. Let ${\hat C}:=\{(0,c):c\in C\}$
consist of the costs augmented with the $0$ slack costs,
and
$${\hat X}\ :=\ \{x\in\Z^{(l+1)\times m\times n}\,,\
\sum_j x^k_{i,j}=s^k_i\,,\ \ \sum_i x^k_{i,j}=d^k_j\,,\
\sum_{k=0}^l x^k_{i,j}=u_{i,j}\,,\ 0\leq x^k_{i,j}\leq u_{i,j}\}\ .$$
It can be shown that the matrix defining the equations of $\hat X$ is an $n$-fold 
product of a fixed bimatrix, hence its Graver basis can be computed in polynomial time, 
see \cite[Chapter 4]{Onn}. The corollary now follows from 
Theorem \ref{Main} applied to $\hat X$ and $\hat C$.
\epr

Next we consider the (integer) {\em three-dimensional transportation problem},
$$\min\{cx\ :\ x\in\Z_+^{l\times m\times n}\ :\ \sum_i x_{i,j,k}=u_{j,k}
\,,\ \sum_j x_{i,j,k}=v_{i,k}\,,\ \sum_k x_{i,j,k}=w_{i,j}\}\ .$$
It is NP-hard even for $l=3$, see \cite{DO}.
But we have the following robustness statement.

\bc{Transportation}
For every fixed $l$ and $m$, the robust three-dimensional transportation problem
$\min_X\max_C$ with box $C$ or $\max_C\min_X$ with list $C$ is polytime solvable.
\ec
\bpr
Note that each variable obeys the bounds
$0\leq x_{i,j,k}\leq\min\{u_{j,k},v_{i,k},w_{i,j}\}$.
It can be shown again that the matrix defining the equations of $X$ is an $n$-fold 
product of a fixed bimatrix, hence its Graver basis can again be computed in polynomial 
time, see \cite[Chapter 5]{Onn}. The corollary then follows again from Theorem \ref{Main}.
\epr

It was shown in \cite{DO} that {\em every} bounded integer program can be
isomorphically represented in polynomial time as some $3\times m\times n$
transportation problems with $l=3$ and some $m$ and $n$. Regarding $m$ as a parameter,
Corollary \ref{Transportation} implies that we can do robust integer programming
in polynomial time for the class of $3\times m\times n$ problems with varying $n$.
Since every integer program belongs to one of these classes for some $m$,
this provides an effective parametrization of robust integer programming.

\end{document}